\documentclass[12pt,leqno]{article}

\title{Integral Positive Ternary Quadratic Forms}
\author{ William C. Jagy  }
\date{}
\begin{document}

\maketitle

\newtheorem{theorem}{Theorem}[section]
\newtheorem{lemma}[theorem]{Lemma}
\newtheorem{cor}[theorem]{Corollary}
\newtheorem{deef}[theorem]{Definition}
\newtheorem{conjecture}[theorem]{Conjecture}

\begin{abstract}
We discuss some  families of integral positive ternary quadratic forms. Our main example is \(f(x,y,z) = x^2 + y^2 + 16 n z^2,\)   where $n$ is positive, squarefree,  and \( n = u^2 + v^2 \) with \( u,v \in \mathbf Z. \)
\end{abstract}

\footnotetext{ADDRESS: Math.Sci.Res.Inst., 17 Gauss Way, Berkeley, CA 94720-5070}   
\footnotetext{EMAIL: jagy@msri.org} 
\footnotetext{SUBJECT CLASSIFICATION: Primary 11E20; Secondary 11D85, 11E12, 11E25}  

\section{Notation}
As in~\cite{bi},~\cite{lehman}, and section 7 of~\cite{glw:thesis}, we let the integer sextuple
\[ \langle a,b,c,r,s,t \rangle  \]
refer to the quadratic form
\[  f(x,y,z) = a x^2 + b y^2 + c z^2 + r y z + s z x + t x y. \]
The Gram matrix for the form is the matrix of second partial derivatives:
\[  
 \left(  \begin{array}{rrr}
  2 a & t & s \\
  t  & 2 b & r \\
  s & r  & 2 c  
\end{array} 
  \right)  .
  \]
So our Gram matrix is symmetric, positive definite, and has integer entries.

We define our  discriminant $ \Delta$ as  half the determinant of the matrix above, so
\[ \Delta = 4 a b c + r s t - a r^2 - b s^2 - c t^2 . \]
   All our forms are positive and primitive ( $\gcd(a,b,c,r,s,t) = 1$). Note that we do allow some of $r,s,t$ to be odd at times. When $r,s,t$ are all even, we refer to the form as classically integral.

\section{Introduction}
In a  1995 letter to J.S.Hsia and R. Schulze-Pillot, Irving  Kaplansky pointed out some simple properties of 
\[ \langle 2,2,4 k^2 + 1, 2,2,0 \rangle  \]
or
\[  f(x,y,z) = 2 x^2 + 2 y^2 + ( 4 k^2 + 1 ) z^2 + 2 y z + 2 z x. \]
When $k$ is odd, then $ f \neq m^2,$ in notation going back to Jones and Pall~\cite{jp}, where this means that all prime factors of $m$ are congruent to $ 1 \pmod 4.$ 

We give the simple proof, while changing the focus to
\[ \langle 2,2,4 n + 1, 2,2,0 \rangle  \]
where $n$ is odd, squarefree, and $n = u^2 + v^2$ in integers. Furthermore the numbers not represented will be all $n m^2.$
\begin{lemma}
Let \( n\) be positive, odd, squarefree, and $n = u^2 + v^2$ in integers. Then
\[ \langle 2,2,4 n + 1, 2,2,0 \rangle \neq n m^2. \]
\end{lemma}
{\bf Proof:}
We have the identity
\[ 2 x^2 + 2 y^2 + ( 4 n + 1 ) z^2 + 2 y z + 2 z x = ( x + y + z)^2 + (x-y)^2 + 4 n z^2. \] That is to say, \( \langle 2,2,4 n + 1, 2,2,0 \rangle\) represents all numbers that can be expressed as \( U^2 + V^2 + 4 n z^2 \) with
\( U + V + z \) even. So, assume we have
\[  U^2 + V^2 + 4 n z^2 = n m^2, \; \; \; U + V + z \equiv 0 \pmod 2. \]
As \(n,m \) are odd, it follows that \( U + V\) is odd, so \(z\) is also odd and nonzero.
Then \[ U^2 + V^2 = n (m^2 - 4 z^2), \] and
 \[ U_1^2 + V_1^2 =  m^2 - 4 z^2 = (m + 2 z)( m - 2 z).\]
Now, \(m + 2 z \equiv 3 \pmod 4,  \; \;m - 2 z  \equiv 3 \pmod 4.\) There is some prime \(q \equiv 3 \pmod 4 \) such that
\( q^{2 i + 1} \parallel m + 2 z.\) However, \((m + 2 z)( m - 2 z) \) is the sum of two squares, 
so we also have \( q^{2 j + 1} \parallel m - 2 z,\) from which it follows that \( q | m,\) a contradiction. \(\bigcirc \)

Our discussion of the genus containing \( \langle 2,2,4 n + 1, 2,2,0 \rangle \) is simplified by
\begin{lemma}
Let \( k\) be any positive  integer. Then \( \langle 1,1,16 k , 0,0,0 \rangle  \)
and \newline \( \langle 2,2,4 k + 1, 2,2,0 \rangle\)   are in the same genus.
\end{lemma}
{\bf Proof:} We use Proposition 4 on page 410 of Lehman~\cite{lehman}, using his terminology and notation, once for each form. Divisor, reciprocal, and level   are defined on page 402, while  conditions we need on the relationship of the form and its reciprocal are given in Proposition 2 on page 403.    

First, we take \(f =  \langle a,b,c,r,s,t \rangle =  \langle 1,1,16 k , 0,0,0 \rangle, \) which has discriminant \(64k,\) level \(64k,\) and divisor \(m=4.\)  Next, we find its reciprocal
\(\phi =  \langle \alpha, \beta,\gamma, \rho, \sigma,\tau \rangle =  \langle 16 k,16k,1 , 0,0,0 \rangle, \) which has discriminant \(1024k^2,\) level \(64k,\) and divisor \(\mu=64k.\) So we have \(a = \gamma = 1.\)

Lehman defines the collection of genus symbols on page 410. As $m=4$ is not divisible by any odd prime or by 16 or 32, none of the genus symbols  $(f| \cdot)$  are defined. As $\mu = 64k$ and $\gamma = 1,$ for any odd prime dividing $k$ we have
$(\phi|p) = (\gamma|p) = (1|p) = 1.$ Then, as $ 16,32 | \mu,$ we have $(\phi|4) = (-1)^{(\gamma - 1)/2} = (-1)^0 = 1,$ then 
$(\phi|8) = (-1)^{(\gamma^2 - 1)/8} = (-1)^0 = 1.$

We need to take a cyclic permutation of variables in our second form to use these results, so, reusing most of the letters, take
\(h =  \langle a,b,c,r,s,t \rangle =  \langle 4k+1,2,2 , 0,2,2 \rangle,\) which has discriminant \(64k,\) level \(64k,\) and divisor \(m=4.\) The reciprocal is  \(\eta =  \langle \alpha, \beta,\gamma, \rho, \sigma,\tau \rangle =  \langle 4, 8 k+1,8k+1,2 , -4,-4 \rangle, \) which has discriminant \(1024k^2,\) level \(64k,\) and divisor \(\mu=64k.\) This time $a = 4k+1$ and $\gamma = 8k+1.$ This works out, insofar as the conditions in Proposition 2 are that $\gcd(a,\gamma) = \gcd(a,m \mu) =\gcd(\gamma, m \mu) = 1.$

Once again, with $m=4,$ Lehman gives no value for any of the genus symbols  $(h| \cdot).$ For any odd prime $p | k,$ we get  $(\eta|p) = (\gamma|p) = (8k+1|p) = (1|p) =  1.$ Then, as $ 16,32 | \mu,$ we have $(\eta|4) = (-1)^{(\gamma - 1)/2} = (-1)^{4k} = 1,$ then 
$(\eta|8) = (-1)^{(\gamma^2 - 1)/8} = (-1)^{8 k^2 + 2 k} = 1.$

We have calculated discriminant, level, and collection of genus symbols for $f,h$ and found agreement, so our two forms are in the same genus by Proposition 4 of ~\cite{lehman}.  \(\bigcirc \)

We introduce a celebrated result of Duke and Schulze-Pillot, which is the Corollary to Theorem 3 in~\cite{dsp}:
\begin{theorem} \label{thm:dsp} Let $q(x_1,x_2,x_3)$ be a positive integral ternary quadratic form. Then every large integer $n$ represented primitively by a form in the spinor genus of $q$ is represented by $q$  itself and the representing vectors are asymptotically uniformly distributed on the ellipsoid $q({\bf x}) = n.$
\end{theorem}

We will also need a short lemma on binary forms:
\begin{lemma}  \label{prim}
 If   all prime factors of a positive integer are  $1 \pmod 4,$ then it can be represented primitively as $x^2 + y^2,$ that is with $\gcd (x,y) = 1.$
\end{lemma}

From Lemma~\ref{prim},  when 
\(n\) is odd, squarefree, and \(n = u^2 + v^2\) in integers, and all prime factors of \(m\) are \(1 \pmod 4\) as well (although \(m\) need not be squarefree), we see that \(n m^2\) is primitively represented by \( \langle 1,1,16 n , 0,0,0 \rangle.\) But Kaplansky's argument has shown that \( \langle 2,2,4 n + 1, 2,2,0 \rangle \neq n m^2. \) It now follows from 
Theorem~\ref{thm:dsp} that  \( \langle 2,2,4 n + 1, 2,2,0 \rangle\) and \( \langle 1,1,16 n , 0,0,0 \rangle,  \) while in the same genus, are in fact in different spinor genera, so there are at least two spinor genera in this genus. 

J. S. Hsia~\cite{hsia} confirmed for the author  that, for both odd and even squarefree \(n=u^2 + v^2, \) the genus of  \( \langle 1,1,16 n , 0,0,0 \rangle\)  has exactly two spinor genera, and that $n$ itself is a {\bf spinor exceptional integer} (a number not represented by one of the spinor genera). He mentioned that the methods were in~\cite{ehh}. 
He also pointed  out his proof that, if there are any spinor exceptions for a genus, there is one that  divides $ 2 \Delta,$ this being Theorem 2 in~\cite{jh2}. Our family shows that the smallest spinor exception can be as large as
\( \Delta / 64.\) 

We return briefly to the base genus, with our \(n=1.\) For  all numbers except odd squares, the number of representations by  \( \langle 1,1,16  , 0,0,0 \rangle  \) is the same as  the number of representations by  \( \langle 2,2,5, 2,2,0 \rangle.\) Then,  for \(k\) {\bf odd},  
\( \; \; r_{\langle 1,1, 16, 0,0,0 \rangle}(k^2) -  r_{\langle 2,2, 5, 2,2,0 \rangle}(k^2) = 4 \; ( -1 | k ) \; k.\)  Complete proofs of these facts have been supplied by Alexander Berkovich~\cite{berkovich} and Wadim Zudilin\cite{zudilin}, in the language of modular forms.
In this situation, the odd squares are called the {\bf splitting integers} for the genus, as the Siegel weighted average representation of the odd squares for one spinor genus disagrees with that of the other spinor genus.
 As it is also possible to calculate the Siegel weighted average of representations for any genus, this allows one to separately calculate
\( r_{\langle 1,1, 16, 0,0,0 \rangle}(j) \) and \(  r_{\langle 2,2, 5, 2,2,0 \rangle}(j) \) for any integer \(j.\) Splitting integers are used in section 2 of~\cite{behh} to correctly partition a genus of ten classes into its spinor genera, five classes each. The characterization of splitting integers as disagreement of representation measures is Corollary 1 on page 3 of~\cite{behh}. An anonymous referee has pointed out that explicit calculation of the difference of representation measures is dealt with in Satz 2 and Korollar 1 of~\cite{rsp:1984dar}.   

\section{A rare phenomenon}
We have mentioned that, with $n$ squarefree and $n=u^2 + v^2,$ the genus of $ \langle 1,1, 16 n, 0,0,0 \rangle  $ has two spinor genera, and $n$ itself is a spinor exception. In this section we prove

\begin{theorem} \label{oddity}
Let \( n\) be positive, odd, squarefree, and $n = u^2 + v^2$ in integers. Then every form in the same spinor genus as
$ \langle 1,1, 16 n, 0,0,0 \rangle  $ also integrally represents $n.$
\end{theorem}
In another section we will prove the same result for even $n.$

The main tool is a {\bf genus-correspondence}, with the first simple properties conjectured by the author, and proved by Wai Kiu Chan~\cite{chan2008}. First we need to describe what we mean by a ternary form representing a multiple of another ternary form.

Suppose we have two positive ternary forms $f,g,$ with Gram matrices $F,G,$ and suppose we have some positive integer $k.$ We will say that $f$ represents $k g$ when there is an integral matrix $P$  such that
\[ P^t \,  F P = k G. \]
The easiest consequence of such a relationship is that, whenever $g$ integrally represents an integer $w,$ it follows that $f$ integrally represents $kw.$
 
Our concern is for the situation when two forms represent prescribed multiples of each other:
 \begin{theorem}{\mbox{\rm (Chan)}} \label{chan} 
Suppose $f_0, g_0$ are positive ternary forms with {\bf integral} discriminant ratio
 $k.$ Suppose that 
$f_0$ represents
 $k g_0$ and 
$g_0$ represents 
$k f_0.$ 
Then, for any 
$f_1 \in \mbox{\rm gen} \; f_0,$ there is at least one  $g_1 \in \mbox{\rm gen} \; g_0$ such that
 $f_1$ represents $k g_1$ and $g_1$ represents $k f_1.$ 
Also, for any $g_2 \in \mbox{\rm gen} \; g_0,$ there is at least one $f_2 \in \mbox{\rm gen} \; f_0$  such that
 $g_2$ represents $k f_2$ and
 $f_2$ represents $k g_2.$  
\end{theorem}
We call this a {\bf genus-correspondence} because it is generally many-to-many, that is, there is generally no well-defined mapping on equivalence classes of forms in either direction.

We are now able to prove Theorem~\ref{oddity}.
Take $n = u^2 + v^2$ to be squarefree and {\bf odd}.
Let $G_0$ be the Gram matrix for   \( g_0 = \langle 1,1,16 n , 0,0,0 \rangle,  \) so that
\[  G_0 \; \; = \; \;  
 \left(  \begin{array}{rrr}
  2 & 0 & 0 \\
  0  & 2 & 0 \\
  0 & 0  & 32 n  
\end{array} 
  \right)  .
  \]
Let $F_0$ be the Gram matrix for   \( f_0 = \langle 1,1,16  , 0,0,0 \rangle,  \) so that
\[  F_0 \; \; = \; \;  
 \left(  \begin{array}{rrr}
  2 & 0 & 0 \\
  0  & 2 & 0 \\
  0 & 0  & 32   
\end{array} 
  \right)  .
  \]

We have \( P^t \, G_0 \, P = n F_0,\) with
\[  P \; \; = \; \;  
 \left(  \begin{array}{rrr}
  u & v & 0 \\
  -v  & u & 0 \\
  0 & 0  & 1  
\end{array} 
  \right)  .
  \]
Note that $\det P = n.$ We take the adjoint $Q$ so that $PQ = QP = nI$ and for that matter 
 $\det Q = n^2.$ We find that 
\( Q^t \, F_0 \, Q = n G_0,\) with
\[  Q \; \; = \; \;  
 \left(  \begin{array}{rrr}
  u & -v & 0 \\
  v  & u & 0 \\
  0 & 0  & n  
\end{array} 
  \right)  .
  \]
Furthermore, the ratio of the discriminants of $f_0,g_0$ is $64n/64=n.$ So $f_0$ represents
 $n g_0$ and 
$g_0$ represents 
$n f_0,$  and Theorem~\ref{chan} applies.

Let $g_1$ be any form in the spinor genus of $g_0,$ written $g_1 \in \mbox{\rm spn} \; g_0.$ According to Lemma~\ref{prim}, for any prime $p \equiv 1 \pmod 4,$ we know that $x^2 + y^2$ and therefore $g_0$ represent $n p^2$ primitively. According to Theorem~\ref{thm:dsp}, when $p$ is  sufficiently large, $n p^2$ is also represented by $g_1.$ From Theorem~\ref{chan}, we know that $g_1$ corresponds with either    \( f_0 = \langle 1,1,16  , 0,0,0 \rangle\) or 
  \( f_1 = \langle 2,2,5,2,2,0 \rangle.\) However, if $f_1$ represented $n g_1,$ then $f_1$ would integrally represent
$n^2 p^2,$ which is a spinor exception for this genus and is not, in fact, represented by $f_1.$ It follows that
$g_1$ represents $n f_0$ and $f_0$ represents $n g_1.$ In particular, $g_1$ integrally represents $n.$ This completes the proof of  Theorem~\ref{oddity}. \(\bigcirc\)

Next, consider any
$g_2 \in \mbox{\rm gen} \; g_0$ but 
$g_2 \notin \mbox{\rm spn} \; g_0.$ Then $g_2$ does not represent  $n,$ as $n$ is a spinor exception for $\mbox{\rm gen} \; g_0.$ So it is not possible for $g_2$ to represent $n f_0.$ From Theorem~\ref{chan}, we find that $g_2$ represents 
$n f_1,$ where   \( f_1 = \langle 2,2,5,2,2,0 \rangle.\) We have chosen to say that this genus-correspondence
respects spinor genus. Formally, we could say this: given a pair of genera with discriminant ratio $k$ and a 
genus-correspondence. Suppose that   $f_3$ represents $k g_3$ and $g_3$ represents $k f_3,$ while 
   $f_4$ represents $k g_4$ and $g_4$ represents $k f_4.$ We say that the genus-correspondence
{\bf respects spinor genus} when $f_3,f_4$ are in the same spinor genus if and only if $g_3,g_4$ are in the same spinor genus. 

We have  extensive numerical support for the following:
\begin{conjecture}\label{conj1}
Given two genera $G_1, G_2$ of positive ternary forms, with integral {\bf squarefree} discriminant ratio and with a genus-correspondence. Suppose that $G_1, G_2$ both have exactly two spinor genera. Then $G_1$ has spinor exceptional integers if and only if  $G_2$ has spinor exceptional integers,  $G_1$ has splitting integers if and only if  $G_2$ has splitting integers, and the genus-correspondence respects spinor genus. When there are spinor exceptions, the regular spinor genera correspond. When there are splitting integers, the spinor genera that have larger (weighted) representation measures for the smallest splitting integers correspond. 
\end{conjecture}

We should emphasize that a genus need not have splitting integers. The best known example is that of 
\( \mbox{\rm gen} \;\langle 1,17,289,0,0,0 \rangle,\) from page 257 of~\cite{bh}. The example with the smallest discriminant (1375) is \( \mbox{\rm gen} \;\langle  1, 5, 70, 5, 0, 0 \rangle,\) just beyond the range of the Brandt and Intrau tables~\cite{bi}. It was rather surprising that splitting integers were not evidently required for a genus-correspondence to respect spinor genus, as there is then no apparent way to label one spinor genus as ``more regular'' than the other.

With less detail and far less evidence, we also offer, for four or more spinor genera,
\begin{conjecture}\label{conj2}
Given two genera $G_1, G_2$ of positive ternary forms, with integral {\bf squarefree} discriminant ratio and with a genus-correspondence. Suppose that $G_1, G_2$  have exactly the same number (some $2^j$) of spinor genera. Then the genus-correspondence respects spinor genus. 
\end{conjecture}

Note that,  with squarefree discriminant ratio and a genus-correspondence, it is still  common for either the genus with larger discriminant or the genus with the smaller discriminant to have fewer spinor genera than the other. Such examples can be quite instructive.

\section{Tornaria's constructions}
Gonzalo Tornaria was kind enough to describe the genus-correspondence, in two situations, as a mapping between forms in some canonical shapes. These mappings do not extend to mappings of equivalence classes. The virtue of this approach is the placing of the genus-correspondence as merely one  variant of Kaplansky's ``descent'' steps, used in preparing~\cite{jks}, and described throughout~\cite{ikclassification}. The similarity to Watson transformations~\cite{glw:trans}  also becomes apparent, although a Watson transformation is a well-defined mapping on equivalence classes of forms, and a Watson transformation does not send a form with some odd prime $p \parallel \Delta$ to a form with $\Delta \neq 0 \pmod p.$ The closest parallel we know involving a Watson transformation is the descent of a form (probably regular) with $\Delta =2592 = 32 \cdot 81$ to one with $\Delta = 32$ that is regular, in that \[ \lambda_9 (  \langle 5,9,17,6,5,3 \rangle  ) =  \langle 1,3,3,1,0,1 \rangle
 .   \]  

We have taken some extra care to show how Tornaria's ascent and descent steps may be viewed as inverses, at least to the extent that they interchange forms in one canonical shape with forms in another canonical shape.

Take an odd prime $p$ and a discriminant such that $\Delta \neq 0 \pmod p.$ Take any form 
\( f_0 = \langle a,b,c,r,s,t \rangle\) with discriminant $\Delta.$ As $f_0$ is isotropic in $\mathbf Q_p,$ we may demand that $c \equiv 0 \pmod p,$
in that such a value is indeed primitively represented by our form. From
\(  \Delta \equiv rst- a r^2 - b s^2 \neq 0 \pmod p    \) we know that $r,s$ are not both divisible by $p.$ If necessary, interchange variables so that $s \neq 0 \pmod p.$ Formally, we have taken the Gram matrix $A_1$ and replaced it by the equivalent $A_2 = P^t \; A_1 P,$ where
\[  P \; \; = \; \;  
 \left(  \begin{array}{rrr}
  0 & 1 & 0 \\
  1  & 0 & 0 \\
  0 & 0  & 1  
\end{array} 
  \right)  .
  \] 
The coefficients become \(  \langle b,a,c,s,r,t \rangle,\) and we simply rename these with the original letters.
So we now have \(  \langle a,b,c,r,s,t \rangle\) with  $c \equiv 0 \pmod p,s \neq 0 \pmod p.$
Next, solve for $k$ in $a + s k \equiv 0 \pmod p,$ then find $A_3 = Q^t \; A_2 Q,$ with
\[  Q \; \; = \; \;  
 \left(  \begin{array}{rrr}
  1 & 0 & 0 \\
  0  & 1 & 0 \\
  k & 0  & 1  
\end{array} 
  \right),
  \]
The new coefficients are \(  \langle a+sk+ck^2,b,c,r,s+2ck,t+rk \rangle.\) Renaming again, we have
\(  \langle a,b,c,r,s,t \rangle\) with  $a,c \equiv 0 \pmod p,s \neq 0 \pmod p,$ this being the first of the two canonical shapes.
 Then we may construct the form
\[ g_0(x,y,z) = \frac{1}{p} \; f_0(px, py, z), \] with coefficients
 \[ g_0 = \left\langle pa, pb, \frac{c}{p},r,s,pt \right\rangle.\]

In the descent direction, let  $\Delta \equiv 0 \pmod p$  and  $\Delta \neq 0 \pmod {p^2},$ or $p \parallel \Delta.$  
Let \( g_1 = \langle a,b,c,r,s,t \rangle\) have discriminant $\Delta.$ This time we need to explicitly require that the form be isotropic in $\mathbf Q_p.$ We then demand that $p^2 | a.$  It follows that
\( \Delta \equiv rst- b s^2 - c t^2 \neq 0 \pmod {p^2}.\)
Thus we know that $s,t$ are not both divisible by $p.$ If necessary, transpose $s,t$ so that $s \neq 0 \pmod p.$ We are taking the Gram matrix $B_1$ and replacing it by $B_2 = P^t \; B_1 P,$ where
\[  P \; \; = \; \;  
 \left(  \begin{array}{rrr}
  1 & 0 & 0 \\
  0  & 0 & 1 \\
  0 & 1  & 0  
\end{array} 
  \right)  .
  \]  
Next, solve for an integer $k$ in $t + s k \equiv 0 \pmod p.$ Construct the matrix
\[  Q \; \; = \; \;  
 \left(  \begin{array}{rrr}
  1 & 0 & 0 \\
  0  & 1 & 0 \\
  0 & k  & 1  
\end{array} 
  \right),
  \]
and take the form with Gram matrix $B_3 = Q^t \; B_2 Q.$
The new coefficients are \(  \langle a,b+rk+ck^2,c,r+2ck,s,t+sk \rangle.\)
 The value $t$ has thus been replaced by $t + s k,$ divisible by $p,$ but without altering the value of $a$ or $s.$ At this point, \( \Delta \equiv- b s^2  \pmod p,\) so that
$p | b.$ We now have our form in the second canonical shape,
\( g_1 = \langle a,b,c,r,s,t \rangle, \) with \(a,b,t\) all divisible by $p,$ indeed  $p^2 | a,$ but $s \neq 0 \pmod p.$ The new form, with discriminant
$\frac{\Delta}{p},$ is given by 
\[ f_1(x,y,z) = \frac{1}{p} \; g_1(x, y, pz), \]
 with coefficients
 \[ f_1 = \left\langle \frac{a}{p}, \frac{b}{p} ,pc,r,s, \frac{t}{p} \right\rangle.\]

\section{Even $n$}
\label{even}
We prove the other case of Theorem~\ref{oddity}, namely 
\begin{theorem} \label{evenity}
Let \( n\) be positive, even, squarefree, and $n = u^2 + v^2$ in integers. Then every form in the same spinor genus as
$ \langle 1,1, 16 n, 0,0,0 \rangle  $ also integrally represents $n.$
\end{theorem}
{\bf Proof:} The genus containing  \( f_0 = \langle 1,1,32 , 0,0,0 \rangle\) consists of three classes, in two spinor genera. The first spinor genus contains the classes  \( \langle 1,1,32 , 0,0,0 \rangle\) and
\( \langle 2,2,9 , 2,2,0 \rangle,\) both of which represent 2. The other spinor genus consists of the single class
 \( \langle 1,4,9 , 4,0,0 \rangle,\) which does not represent 2 or any $2 m^2.$

With $n$ even,  $ g_0 = \langle 1,1, 16 n, 0,0,0 \rangle  $ represents
 \( \frac{n}{2} \cdot \langle 1,1,32, 0,0,0 \rangle,\) so that   \( \langle 1,1,32, 0,0,0 \rangle\) also represents
 \( \frac{n}{2} \cdot \langle 1,1,16 n, 0,0,0 \rangle,\) and there is thus a genus-correspondence. Consider some $g_1 \in \mbox{\rm spn} \; g_0.$ From Lemma~\ref{prim}, for any prime $p \equiv 1 \pmod 4,$ we know that $x^2 + y^2$  represents $(n/2) p^2$ primitively, denote this $(n/2) p^2 = a^2 + b^2, \; \; \gcd(a,b) = 1.$ As $a^2 + b^2$ is odd, it follows that 
$ \gcd(a-b,a+b) = 1$ as well. So we have the primitive representation
$(a-b)^2 + (a+b)^2 = n p^2,$ which tells us that $g_0$ primitively represents $n p^2.$ When $p$ is sufficiently large, 
Theorem~\ref{thm:dsp} tells us that $g_1$ represents $n p^2.$ By Theorem~\ref{chan}, we know that $g_1$ corresponds with at least one of the three forms in the genus of $f_0.$ However, if   \( \langle 1,4,9 , 4,0,0 \rangle\) should represent $\frac{n}{2} g_1,$ it would follow that \( \langle 1,4,9 , 4,0,0 \rangle\) represented the integer $\frac{n^2 p^2}{2},$ which is of the form $2 m^2.$ It follows that $g_0$ represents either  \( \frac{n}{2} \cdot \langle 1,1,32, 0,0,0 \rangle\) or
 \( \frac{n}{2} \cdot \langle 2,2,9,2,2,0 \rangle.\) In either case $g_0$ represents the integer $n.$ \(\bigcirc\)

We pause to discuss the influence of Conjecture~\ref{conj1}. It was necessary to have a separate proof for even $n$ because $4 m^2$ is not a spinor exception for the genus containing  \( \langle 1,1,16, 0,0,0 \rangle.\) If we had known a proof of Conjecture~\ref{conj1}, we could simply have said that any form in the same spinor genus as   \( \langle 1,1,16n, 0,0,0 \rangle\) represents  \( n \cdot \langle 1,1,16 n, 0,0,0 \rangle.\) Similarly, we would not have needed any invocation of Theorem~\ref{thm:dsp}, which can become unusable if primitive representations of desirable numbers are not available.

Conjecture~\ref{conj1} would be an even bigger help in the following related pair of examples, where the conjectured behavior has simply not been proved, although checked as correct for $n \leq 200.$ One situation is $n= u^2 + u v + 4 v^2$ squarefree, with the genus of \( \langle 1,4,225 n , 0,0,1 \rangle.\)  Second,  $n= 2 u^2 + u v + 2 v^2$ squarefree, and the genus of \( \langle 2,2,225 n , 0,0,1 \rangle.\) In these cases $n$ is allowed odd or even. The  ``base'' genus has  four forms in two spinor genera of two classes each:\( \langle 1,4,225  , 0,0,1 \rangle\)
and \( \langle 1,15,60  , 15,0,0 \rangle\) are in one spinor genus, \( \langle  6, 6, 25, 0, 0, 3  \rangle\) and
\( \langle  9, 10, 10, 5, 0, 0  \rangle\) are in the other. The spinor exceptions are of the form $\mu^2,$ where all prime factors of $\mu$ are $1,2,4,8 \pmod{15},$ and 2 itself is included. As $9 \mu^2$ and $25\mu^2$ are not spinor exceptions, to deal with $n$ divisible by $3,5,15,$ we would first need to calculate the genera of \( \langle 2,2,675 , 0,0,1 \rangle,\) \( \langle 2,2, 1125 , 0,0,1 \rangle,\) and \( \langle 1,4,3375 , 0,0,1 \rangle.\)

 

\section{Involutions} 
We return to odd squarefree $n = u^2 + v^2$ and the genus of  \( \langle 1,1,16 n , 0,0,0 \rangle.\) As long as $n \leq 505,$ a few interesting things happen. First, the two spinor genera in the genus have the same number of equivalence classes of forms. Second, for each class $f,$ there is a single class $g$ with $f \neq g,$ such that
$f$ represents $4g$ and $g$ represents $4f,$ while $g$ never lies in the same spinor genus as $f.$ So ``involution'' seems a good term for this, as we have a bijection that interchanges the two spinor genera.

A similar thing happens in these two situations, from the last paragraph of section~\ref{even}: first, $n= u^2 + u v + 4 v^2,$ with \( \langle 1,4,225 n , 0,0,1 \rangle,\) or second,  $n= 2 u^2 + u v + 2 v^2,$ with \( \langle 2,2,225 n , 0,0,1 \rangle,\) while we keep $n$ squarefree, but add the restriction that $n$ not be divisible by 3 or 5. There are indeed two spinor genera, and they are the same size, checked for $n \leq 200.$ The worthwhile detail is that we get one involution where each $f$ has a single $g \neq f$ such that $f$ represents $9g$ and $g$ represents $9f,$ so that is one involution. But there is a different involution where  $f$ represents $25g$ and $g$ represents $25f.$ Both 9 and 25 interchange spinor genera. Nothing special occurs with 4. 

This last conjecture has not been checked as thoroughly, but is worthwhile for suggesting possibilities with four spinor genera. In~\cite{bh}, there is a genus with four spinor genera described, containing the form called
\( B^1 = \langle 1,20,400 , 0,0,0 \rangle.\) The spinor genera all have three classes. There are two families of spinor exceptions, $5m^2,$  all prime factors of $m$ being $1 \pmod 4,$ and $\phi^2,$ where all prime factors  of $\phi$ are  $1,3,7,9\pmod {20}.$

This first step has been checked for $n \leq 1189=29 \cdot 41.$ Let $n$ be squarefree, and all prime factors of $n$ be either $1 \pmod {20}$ or $9 \pmod {20}.$ Then the genus of \( \langle 1,20,400 n , 0,0,0 \rangle\) has four spinor genera of equal size. Either $n=u^2 + 20 v^2$ or $n=4u^2 + 5 v^2,$ and it is easy to check that \( \langle 1,20,400 n , 0,0,0 \rangle\) represents either \( n \cdot \langle 1,20,400, 0,0,0 \rangle\) or \( n \cdot \langle 4,5,400, 0,0,0 \rangle.\) In turn, the relevant form in the ``base'' genus   represents \( n \cdot \langle 1,20,400n, 0,0,0 \rangle.\) This extends to a genus-correspondence. With all as described, this genus-correspondence respects spinor genus.

Let us label the four spinor genera. Let $A_n$ be regular, let $B_n \neq 5nm^2,$ let $C_n \neq n\phi^2,$ finally  $D_n \neq 5nm^2,n\phi^2.$ These next items have been checked only as far as $n \leq 61.$ There are involutions with multiplier 25, these interchange $A_n$ with $C_n,$ and then interchange $B_n$ with $D_n.$ 

In comparison, with multiplier 4, any form in $A_n$ corresponds with a single one in $D_n,$ but with two forms each in $B_n,C_n.$ Similar comments apply beginning with any of the four spinor genera. So multiplier 4 does give an identifiable involution, ($B_n$ matches with $C_n,$) but the behavior is not as clean as that with multiplier 25.

Finally, we explain the restriction on $n$ itself. If $n$ is a number that is represented by both the binary forms 
$x^2 + 20 y^2$  and $4 x^2 + 5 y^2,$ such as $n=21,$ then   \( \langle 1,20,400 n , 0,0,0 \rangle\) represents {\bf both} \( n \cdot \langle 1,20,400, 0,0,0 \rangle\) {\bf and} \( n \cdot \langle 4,5,400, 0,0,0 \rangle,\) so that it is not possible to have a genus-correspondence that respects spinor genus, even if the resulting genus does actually possess four spinor genera.


\section{Acknowledgements}
This article is a selection from the author's talk during a Focused Week at the University of Florida, Gainesville, in March 2010. The author would like to thank the organizers for a delightful event.

We would also like to thank A. Berkovich, W. K. Chan,  J. S. Hsia, G. Tornaria, and W. Zudilin for fundamental contributions to this article and permission to include these. We especially wish to thank W. Zudilin for his careful reading and extensive comments.

\bibliography{terntricks}

\bibliographystyle{plain}

\end{document}